\documentclass[hidelinks,a4paper]{article}
\usepackage[T1]{fontenc}
\usepackage[utf8]{inputenc}
\usepackage{lmodern}
\usepackage{microtype}

\makeatletter

\usepackage[total={5.8in, 7.8in}, asymmetric, bindingoffset=0.4in]{geometry}
\usepackage{amsmath}
\usepackage{amssymb}
\usepackage{amsthm}
\usepackage{mathtools}
\usepackage{enumerate}
\usepackage[dvipsnames]{xcolor}

\usepackage{verbatim}

\usepackage{multirow}
\usepackage{empheq}
\usepackage{listings}
\usepackage{color}

\usepackage[square,comma,numbers]{natbib}
\usepackage{hyperref}
\usepackage{graphicx} % Required for including pictures
\usepackage[font=small,labelfont=bf]{caption}
\numberwithin{equation}{section}
\usepackage{pgfplots}
\pgfplotsset{width=10cm,compat=1.9}

\usepackage[symbol]{footmisc}

\usepackage{fancyhdr}
\date{\vspace{-1em}\normalsize{\today}}

\def\cA{{\mathcal A}}

\def\cC{{\mathcal C}}

\def\cH{{\mathcal H}}
\def\cI{{\mathcal I}}
\def\cJ{{\mathcal J}}
\def\cK{{\mathcal K}}
\def\cL{{\mathcal L}}
\def\cM{{\mathcal M}}

\def\cO{{\mathcal O}}
\def\cP{{\mathcal P}}

\def\cX{{\mathcal X}}

\def\E{\mathbb{E}}

\def\R{\mathbb{R}}

\def\L{\mathbb{L}}

\theoremstyle{plain}
\newtheorem{theorem}{Theorem}[section]
\newtheorem{lemma}[theorem]{Lemma}
\newtheorem{corollary}[theorem]{Corollary}
\newtheorem{proposition}[theorem]{Proposition}
\newtheorem{assumption}[theorem]{Assumption}
\newtheorem{definition}[theorem]{Definition}
\newtheorem{example}[theorem]{Example}
\newtheorem{remark}[theorem]{Remark}

\def\be{\begin{equation}}
\def\ee{\end{equation}}
\def\ba{\begin{align}}
\def\ea{\end{align}}
\def\ba*{\begin{align*}}
\def\ea*{\end{align*}}

\def\fF{{\mathfrak F}}

\def\d{\mathrm{d}}
\def\eps{\epsilon}

\def\om{\omega}
\def\pmu{\partial_\mu}
\def\bu {\bar{u} }

\def\vt{\vartheta}

\def\muk{\mu_k}
\def\nuk{\nu_k}

\def\mud{\mu_\delta}

\def\Fed{\Phi_{\eps,\delta}}
\def\fed{\psi_{\eps,\delta}}
\def\ped{\phi_{\eps,\delta}}
\def\mued{\mu_{\eps,\delta}}
\def\nued{\nu_{\eps,\delta}}
\def\ked{\kappa_{\eps,\delta}}
\def\etaed{\eta_{\eps,\delta}}

\def\td{t_\delta}

\def\taued{\tau_{\eps,\delta}}
\def\ted{t_{\eps,\delta}}

\def\sd{s_\delta}

\def\sed{s_{\eps,\delta}}
\def\oo{\overline{\cO}}
\def\Ied{\cI_{\eps,\delta}}
\def\Jed{\cJ_{\eps,\delta}}

\def\hc{\hat{c}}

\title{Viscosity Solutions of the Eikonal Equation
on the Wasserstein Space\footnote{Partially supported 
by the National Science Foundation grant
 DMS 2106462. Authors thank  
 Professors Martin Larsson, Jianjun Zhou and Jianfeng Zhang for valuable comments,
 and to Professor Samuel Daudin for pointing 
 out an error in an earlier  version.}}

\author{H. Mete Soner\footnote{Department of Operations Research and Financial
Engineering, Princeton University, Princeton, NJ, 08540, USA, email: 
{\tt soner@princeton.edu}. }
\and Qinxin Yan\footnote{Program in Applied and Computational
Mathematics, Princeton University, Princeton, NJ, 08540, USA, email: 
{\tt qy3953@princeton.edu}. }}

\date{\today}

\begin{document}
\maketitle
\abstract{Dynamic programming equations for 
mean field control problems with a separable 
structure are Eikonal type equations on the Wasserstein space.  
Standard differentiation using linear derivatives
yield a direct extension of the classical viscosity
theory.  We use Fourier representation of the 
Sobolev norms on the space of measures,
together with the standard techniques from
the finite dimensional theory to prove
a comparison result among semi-continuous
sub and super solutions, obtaining a unique
characterization of the value function. }
\vspace{3pt}

\vspace{3pt}


\section{Introduction}
\label{sec:intro}

We consider the Hamilton-Jacobi equations related 
to mean-field control problems in which the state 
process $X_t$ taking values in the 
$d$-dimensional Euclidean space $\R^d$
has the following simple dynamic structure,
$$
\d X_u = \alpha_u \, \d u + 
\sigma \, 
\d W_u,
$$
where $\alpha$ is the \emph{control process}
adapted to the information flow but unrestricted otherwise, 
positive square matrix $\sigma$ is the diffusion coefficient,
and $W_t$ is a standard Brownian motion.
The cost functional of these problems have a separable structure given by,
$$
J(\alpha):= \int_t^T [\ell(u,\cL(X_u)) + \frac12\, \E |\alpha_u|^2 ]\, \d u
+ g(\cL(X_T))\, ,
$$
where $\ell, g$, are given functions, and $\cL(X_u) \in \cP(\R^d)$ is the law of the random
variable $X_u$. Let $v(t,\mu)$ be the value function
defined by,
$$
v(t,\mu):= \inf_{\alpha_\cdot} \, J(\alpha),\qquad
\cL(X_t)=\mu.
$$
By appropriately scaling time and space, we assume that $\sigma$
is the identity matrix. Then, the corresponding dynamic programming equation
is given by,
\be
\label{e.hj}
-\partial_t v(t,\mu) + H(\mu, \partial_\mu v(t,\mu)) = \ell(t,\mu),
\ee
where the function $\partial_\mu v(t,\mu)(\cdot)$ 
is the linear derivative of $v$ with respect to $\mu$
as defined in Section \ref{sec:notations} below, and 
for a twice differentiable function $\kappa$ and
a probability measure $\mu$,
\be
\label{e.H}
H(\mu,\kappa)= -\frac12 \mu(\Delta \kappa) 
+\frac12  \mu( |\nabla \kappa|^2),
\ee
and $\mu(f)$ is the action of the measure $\mu$
on the function $f$.  

Under natural assumptions on $\ell, g$
(cf.~Assumption \ref{a.main}, below),
dynamic programming
holds and the value function is a viscosity solution of \eqref{e.hj}
using the standard notion of linear derivative.
Many similar results of this type
have already been proved in far greater generality.
We refer the reader to our previous paper \cite{SY}
for these types of results, and the relevant references
therein.

Mean-field optimal control
problems  are
part of the exciting general program of
Lasry \& Lions \cite{LL1,LL2,LL} as outlined by Lions
during his College de France lectures \cite{Lions}.
Similar type of differential games were
also independently introduced by Huang,
Malham\'{e}, \& Caines~\cite{HMRC}, and 
we refer the reader to
the classical
book of Carmona \& Delarue \cite{CD},
to the lecture notes of Cardaliaguet \cite{C}
for detailed information and more references.

Our central goal is the characterization
of the value function as the unique weak solution of  \eqref{e.hj}.
While the impressive paper of Cardaliaguet \emph{et.~al.}~\cite{CDLL}
provides regularity results for  mean field games,
it is well known that dynamic programming equations
in general do not admit classical solutions, and 
we naturally consider the celebrated viscosity solutions of Crandall \& Lions \cite{CL,CEL, CIL, FS}.  
However, in infinite dimensions the  Hamiltonian
is often not defined when the derivative of the solution is  not in the domain of 
corresponding unbounded operators, as explained in  
the excellent book of Gozzi \& Swiech \cite{GS}.  Thus, the original definition must
be modified, and there are several alternatives.
Among those we pursue the standard definition of a viscosity solution using the linear derivative on the convex set of probability measures, as we have done in our earlier paper \cite{SY}. 

Our main contribution Theorem \ref{t.compare} is  a comparison result for the dynamic programming equation  \eqref{e.hj} among all  semi-continous sub and supersolutions.  
More general results in this direction has already
been proved by Cosso {\emph{et.~al.}}~\cite{CGKPR},
and more recently by Daudin, Seeger \cite{DS} and by Daudin, Jackson \& Seeger \cite{DJS}.
However, we use a different and an alternate
technique developed 
in \cite{SY} based on negative Sobolev norms and their Fourier representations, 
but without using the strong structure imposed on the controls in \cite{SY}.
An important ingredient is the Lipschitz regularity
in the negative Sobolev norms of the value
of optimal control problems
with  smooth coefficients proved in
Proposition \ref{pr:samuel}.  These estimates
were first  used in \cite{SY} in this context.  In the separable 
structure that we consider, it is proved more generally
by Daudin, Delarue \& Jackson \cite{DDJ} using 
the theory of elliptic equations, and were then used in \cite{DJS} to obtain a general comparison
on the $d$-dimensional torus.
We also leverage this Lipschitz regularity of the value functions
and the techniques of  \cite{SY} 
to prove the general comparison result  Theorem \ref{t.compare}
on the whole $\R^d$,
under a weak
regularity condition Assumption \ref{a.main}.

Properties of the solutions of
Hamilton-Jacobi equations on the spaces  of probability measures
have been actively researched in the past two decades.
A milestone in the these studies
is the  lifting introduced by Lions in \cite{Lions}.
This approach maps the problem to an $\L^2$ 
space and connects to the earlier results exploiting the
Hilbert structure, and is further developed
in several papers including \cite{BCFP,BCP,CGKPR1,PW1,PW}.
Additionally, the novel \emph{Lions derivative} and its properties 
are explored in the book of Carmona \& Delarue \cite{CD}.  

As mentioned earlier, \cite{CGKPR} proves a very general
comparison result by extending the deep techniques developed
by Lions  \cite{lions1983}
to the Wasserstein space and
covering essentially all convex Hamiltonians.
Two recent papers \cite{DS,DJS} also prove comparison results
with techniques closer to ours.
While an intriguing new definition
together with the differentiable structure of the Wasserstein two
metric is used in \cite{DS},
\cite{DJS} uses amalgam of deep techniques
including the negative Sobolev norms
and a change of variables introduced
in \cite{BEZ1} to prove several interesting results
on the $d$-dimensional torus.
Also a general Crandall-Ishi type result
is proved in \cite{BEZ1} using the negative
Sobolev norms introduced in \cite{SY} and in this paper.
Additionally, in another recent study \cite{Be}
related to stochastic optimal transport, Bertucci
introduces a highly original new definition of viscosity solutions and proves general comparison principles.
An interesting approach developed by Gangbo \& Swiech \cite{GS} 
and Marigonda \& Quincampoix \cite{MQ}, and Jimenez \emph{et.al.}~\cite{JMQ}
utilizes deep connections to geometry.
Gangbo \& Tudorascu \cite{GT} 
connects this method to Lions’ lifting. 
Cecchin and Delarue \cite{CD} uses Fourier approximations of the 
measures and exploits the semi-concavity, and provides 
an excellent overview of the problem.
In our earlier work \cite{MVJump,SY},
we have used the direct definition of the 
viscosity solutions and 
employed the classical techniques.

Alternatively, projections of these equations
to finite-dimensional
spaces yield approximate equations that
can be directly analyzed by classical results \cite{CIL}.
A second-order problem studied in \cite{CKLS}
provides a clear example of this approach
as its projections
exactly solve the projected finite dimensional equations.
However, in general these projections are only approximate solutions,
and clearly one has to 
effectively control the approximation error
to obtain relevant results.
This is achieved by Cosso \emph{et.al.}~\cite{CGKPR} 
via the  smooth variational
principle together with  
Gaussian smoothed Wasserstein metrics.
Bayraktar \emph{et.al.}~\cite{BEZ} use a different
approach, and Gangbo 
\emph{et.al.}~\cite{GMS} studies the pure projection problem. 

Other highly relevant studies include Wu \& Zhang \cite{WZ}
for path-dependent equations,
Conforti \emph{et.al.}~\cite{CKT} for gradient flows, 
 and Talbi \emph{et.al.}~\ \cite{TTZ1,TTZ2} for mean-field stopping problems.
Additionally,  Ambrosio  \& Feng \cite{AF}, 
and Feng \& Katsoulakis \cite{FK} study 
the closely connected Hamilton Jacobi equations on metric spaces.

The paper is organized as follow.  
General  structure and notations are given in the next section.
In Section~\ref{sec:problem} we briefly define the problem, and
state the standing assumption. Viscosity solutions are defined in Section~\ref{sec:vis},
and the main comparison result Theorem \ref{t.compare} is  
stated and proved in
Section~\ref{sec:main}. In the Appendices, we prove 
a technical lemma and outline the
proof of the
regularity result proved in \cite{DDJ}.
\section{Notations}
\label{sec:notations}

In this section, we summarize the notations and 
known results  used in the sequel.
We denote the dimension of the ambient space by $d$,
and  the finite horizon by $T>0$.
$\cM(\R^d)$ is the set of all Radon measures,
$\cP(\R^d)$ is the set of
probability measures on $\R^d$, and 
$$ 
\cP_2(\R^d) := \{\ \mu \in \cP(\R^d)\ :\
 \small{\int} |x|^2 \ \mu(\d x)< \infty\ \}.
$$
We write $\cM, \cP, \cP_2$ when
the ambient space is clear or irrelevant.
We endow all these spaces with the
the weak* topology
and write $\muk\rightharpoonup \mu$ 
when $\muk$ weak$^*$ converges to $\mu$.

We set $\cO:=(0,T)\times\cP_2$
and endow $\oo:= [0,T] \times \cP_2$ with the product of
Euclidean and weak* topologies.
We utilize the local compactness of $\oo$.
Indeed, set
\be
\label{eq:q}
\vt(\mu):= \mu(q)=\int q(x)\ \mu(\d x), \quad \mu \in \cP_2,
\qquad 
q(x):= \sqrt{1+|x|^2\, }, \quad x \in \R^d.
\ee
Then,  for any constant $c>0$,
the sublevel set  $\{ (t,\mu) \in \oo\ :\
\vt(\mu) \le c\ \}$ is compact.

For metric spaces $E,F$, $\cC(E\mapsto F)$ denotes
the $F$-valued continuous functions on $E$.  We write
$\cC(E)$ when $F=\R$ and  $\cC_b(E)$ for the bounded ones.
For a positive integer $n$, $\cC^n(\R^d)$ is the set
of all $n$-times continuously differentiable, real-valued
functions, and we set
$$
\cC_*:= \cC_*(\R^d) = \{\ f \in \cC(\R^d)\ :\
|f(x)| \le c(1+|x|^2), \ \text{for some constant}\ c\}.
$$
It is clear that $\int f \d \mu$ is well-defined
for the pair $\mu \in \cP_2$, $f \in \cC_*$, and
whenever defined we write
$\mu(f)$ for the integral $\int_{\R^d} f(x)\mu(\d x)$.
We also use the notation,
\be
\label{eq:c2star}
\cC^2_*:= \{\ f \in \cC^2(\R^d) \ :\
f, \ |\nabla f|^2 \in \cC_*, \  D^2f \in \cC_b \ \}.
\end{equation}

Using the standard notion of linear derivative
on the convex set $\cP_2$, we say that
$\varphi \in \cC(\cP_2)$ is 
\emph{continuously differentiable}
if there exists 
$\pmu \varphi \in \cC(\cP_2 \mapsto \cC_*)$
satisfying,
$$
\varphi(\nu) =
\varphi(\mu) + \int_0^1\, (\nu-\mu)(\pmu \varphi(\mu + \tau (\nu-\mu))\, 
\d \tau,\qquad
\forall\ \mu,\nu \in \cP_2.
$$
Clearly,  $\partial_\mu \varphi(\mu)\in \cC_*$ has many representatives. 
However, when $\pmu \varphi(\mu)$ is twice differentiable, 
then $\mu(\Delta \pmu \varphi(\mu))$, and 
$\mu(h(\nabla \pmu \varphi(\mu)))$ with
any continuous function $h$ and appropriate integrability are independent of this choice,
see for instance \cite{CKLS}[Appendix B].
For $\psi \in \cC(\oo)$ and $(t,\mu)\in \cO$, 
$\partial_t\psi(t,\mu)\in \R$ is the time derivative
evaluated at $(t,\mu)$, and  $\partial_\mu\psi(t,\mu) \in \cC_*$
is the derivative in the $\mu$-variable.

We consider the Fourier basis given by,
$$
e(x,\xi):= (2\pi)^{-\frac{d}{2}}\ e^{i\xi\cdot x},\qquad
x \in \R^d,\ \xi\in \R^d,
$$
where $i =\sqrt{-1}$ and $z^*$ is the complex conjugate of $z$.
Then,  for any $f \in \L^2(\R^d)$,
$$
f(x)= \int_{\R^d}\ \fF(f)(\xi) e(x,\xi)\ \d \xi,
\quad
\text{where}
\quad
\fF(f)(\xi):= \int_{\R^d}\ f(x) e^*(x,\xi)\ \d x, \ \
x,\xi \in \R^d.
$$ 
For $ s \in \R$, $\cH_s(\R^d)$ is the classical Sobolev space with 
fractional derivatives \cite{AF,Sob}.  Then,
$$
\|f\|^2_s:= \|f\|^2_{\cH_s(\R^d)}=
 \int_{\R^d}\ (1+|\xi|^2)^s\, |\fF(f)(\xi)|^2\ \d \xi.
$$
Moreover, for $s >k+\frac{d}{2}$, 
$\cH_s(\R^d)$ continuously embeds into 
$\cC_b^k(\R^d)$. Therefore, for $s>\frac{d}{2}$,
$\cM(\R^d) \subset \cH_{-s}(\R^d)$, and $\|\cdot\|_{-s}$ 
is well defined on $\cM(\R^d)$.   Then, for $\eta \in \cM(\R^d)$,
$$
\|\eta\|^2_{-s}= \int_{\R^d} (1+|\xi|^2)^{-s} |\fF(\eta)(\xi)|^2\ \d \xi,
\ \ \text{where}\ \ 
\fF(\eta)(\xi)= \int_{\R^d}\ e^*(x,\xi)\ \eta(\d x), \ \
\xi \in \R^d.
$$
Moreover, by duality,
$$
\|\eta\|_{-s} = \sup\{\  \eta(\psi)\ :\ \psi \in \ \cH_s(\R^d), \,
\|\psi\|_s \le 1\, \}.
$$

We use the choice
\begin{equation}
\label{eq:star}
n_*:=n_*(d)=3+ \lfloor\frac{d}{2}\rfloor,\qquad
\varrho:=\|\cdot\|_{-n_*},
\end{equation}
where $ \lfloor a \rfloor$
is the integer part of a real number $a$.
As $n_* > 2 + \frac{d}{2}$, $\cH_{n_*}(\R^d) \subset \cC^2_b(\R^d)$, 
and by Morrey's inequality
there is a constant $k_d$ depending only on the
dimension such that (see for instance, \cite{A}[Chapter 4])
\begin{equation}
\label{eq:sobolev}
\| \kappa\|_{\cC^2(\R^d)} \le k_d  \| \kappa\|_{n^*},
\qquad \forall\ \kappa \in \cH_{n^*}(\R^d).
\end{equation}

\section{McKean-Vlasov control }
\label{sec:problem}

Let $v(t,\mu)$ be the value function of
the \emph{McKean-Vlasov optimal control} problem defined in the 
Introduction by using all square integrable, adapted  controls. 
For more information, we refer the reader to Chapter 6 in \cite{CD} and \cite{Dau,SY}.  
In particular, the recent paper of Daudin \cite{Dau}
outlines the connections between several formulations 
and prove existence of optimal feedback controls.

Following is the only assumption of the paper.

\begin{assumption}
\label{a.main} We assume that $\ell :\oo \mapsto \R$
is bounded and continuous in the product of Euclidean and 
weak$^*$ topologies, and $g :\cP_2 \mapsto \R$
is bounded and weak$^*$ continuous.
We additionally assume that, there exists a sequence of
smooth functions 
$(\ell_n,g_n)$ approximating $(\ell,g)$
uniformly, a constant $k_*>0$, a modulus $\om$ {\rm{(}}i.e.,
$\om:\R_+ \mapsto \R_+$ is a continuous function with $\om(0)=0${\rm{)}}, 
and  constants $c_n$, such that
for each  $n$, $t,s \in [0,T]$, and  $\mu \in \cP_2$,
$$
|\ell_n(t,\mu)| +|g_n(\mu)|  \le k_*,
\quad
|\ell_n(t,\mu) - \ell_n(s,\mu)| \le k_*\, \om(|t-s|),
$$
\begin{equation*}
\label{eq:elln}
\|\partial_\mu \ell_n(t,\mu)\|_{\cH_{2n^*}(\R^d)}+
\|\partial_\mu \ell_n(t,\mu)\|_{\cC^{2n^*}(\R^d)}+
\|\partial_\mu g_n(\mu)\|_{\cH_{2n^*}(\R^d)}+
\|\partial_\mu g_n(\mu)\|_{\cC^{2n^*}(\R^d)}
\le c_n.
\end{equation*}
\end{assumption}

Above assumption is satisfied by a large class
of functions, and the choice $2n^*$
is arbitrary but does not decrease the
generality. Below we provide a natural class
of such functions.
In fact, regularization techniques 
developed in \cite{CD} can be used to construct the approximating sequence
directly under assumptions on $(\ell,g)$.

\begin{example}
\label{ex:ell}
{\rm{Consider a function $\ell(\mu)=L(\mu(f))$
for some $L\in \cC_b(\R), f\in \cC_b(\R^d)$.
Additionally, assume that $L$ is Lipschitz,
and $f$ is square integrable.
Then, by mollification  one can construct
smooth functions $(L_n,f_n)$
approximating $(L,f)$ uniformly, and
satisfying $\|f_n\|_{\cH_{2n^*}(\R^d)}+\|f_n\|_{\cC^{2n^*}(\R^d)} \le c_n$,
$$
\sup_n ( \| L_n \|_{\cC^1}  +\|f_n\|_\infty) \le \|L\|_\infty
+\|L^\prime\|_\infty+ \|f\|_\infty =:k_*.
$$ 
Moreover, as $\partial_\mu \ell_n(t,\mu)(x) =
L_n^\prime(\mu(f_n)) f_n(x)$ for $x \in \R^d$,
$$
\|\partial_\mu \ell_n(t,\mu)\|_{\cH_{2n^*}(\R^d)}
\le k_*\ \| f_n\|_{\cH_{2n^*}(\R^d)},
\qquad
\|\partial_\mu \ell_n(t,\mu)\|_{\cC^{2n^*}(\R^d)}
\le k_*\ \| f_n\|_{\cC^{2n^*}(\R^d)}.
$$
Thus, $\ell(\mu)$ satisfies the
above assumptions. More generally, 
a natural class of functions for the
above assumption is given by
$\ell(t,\mu)= L(t,\mu(f_1(t,\cdot)),\ldots,\mu(f_m(t,\cdot)))$
for some functions $L, f_1,\ldots,f_m$
satisfying appropriate conditions. }} 
\end{example}

Let $(\ell_n,g_n)$ be as in the Assumption \ref{a.main},
and $v_n$ be the value function of the optimal control
problem with running cost $\ell_n$ and terminal cost $g_n$,
and same dynamics as in the original problem.
The following regularity of  $v_n$ is
essentially proved in \cite{DDJ}[Proposition 3.2] improving
a similar result proved in \cite{SY}[Theorem 4.2].

\begin{proposition}[Proposition 3.2 \cite{DDJ}]
\label{pr:samuel} 
Let $\varrho$ be as in \eqref{eq:star}.
Under the Assumption \ref{a.main},
there exists constants $\hat{c}_n$ such that
\begin{equation}
\label{eq:lip}
|v_n(t,\mu)-v_n(t,\nu)| \le \hat{c}_n\ \varrho(\mu-\nu),
\qquad \forall \ t \in [0,T],\ \mu,\nu \in \cP_2.
\end{equation}
\end{proposition}

Proposition 3.2 in \cite{DDJ} proves exactly
the above estimate but in the $d$-dimensional 
torus. However, their proof can be directly adopted 
to the current context with no changes.  
As the above estimate is used
centrally in our proofs, for the convenience
of the readers we provide an
outline proof of the above result in the Appendix.

\begin{corollary}
\label{cor:continuous}
Under the Assumption \ref{a.main},
$v_n, v \in \cC_b(\oo)$, i.e.,
both $v_n$ and $v$ are bounded and are continuous
in the product of Euclidean and weak$^*$ 
topologies.
\end{corollary}
\begin{proof}
The continuity of $v_n$  in the time variable
is straightforward \cite{SY}. The above Lipschitz continuity in $\varrho$
and Lemma \ref{lem:equiv} implies that $v_n \in \cC_b(\oo)$.
The uniform convergence of $(\ell_n,g_n)$ to $(\ell,g)$
implies that $v_n$ converges to $v$ uniformly and
therefore  $v\in \cC_b(\oo)$ as well.
\end{proof}

\section{Viscosity Solutions}
\label{sec:vis}

We start by defining the class of
test functions used in the definition of
the viscosity solutions.
\begin{definition}
\label{def:smooth} {\rm{A continuous function
$\varphi \in \cC(\oo)$ is called a}}
test function {\rm{if there exists a 
version of $\pmu \varphi$ such that
the map 
$$
(t,\mu,x) \in \oo \times \R^d\mapsto \pmu \psi(t,\mu)(x)
$$
is continuous, and 
$\pmu \varphi(t,\mu)\in \cC^2_*$
for every $(t,\mu) \in \cO$.
Let $\cC_s(\cO)$ be 
the set of all smooth test functions.
}}\end{definition}

We can now directly define the notion
of viscosity solutions \cite{CL,CEL,CIL,FS}.  Recall that
we endow $\oo$ with the product of Euclidian
and weak$^*$ topologies. 
 \begin{definition}
 \label{def:viscosity}
 We say that an upper semicontinuous
 function $u :\oo \mapsto \R$ is a
{\rm{viscosity subsolution}}of \eqref{e.hj} if for
 every test function $\varphi \in \cC_s(\oo)$
 we have
$$
 -\partial_t \varphi(t_0,\mu_0) 
 +H(\mu_0,\partial_\mu \varphi(t_0,\mu_0)) \le \ell(t_0,\mu_0),
$$
 at every $(t_0,\mu_0) \in \cO$ satisfying
 $(u-\varphi)(t_0,\mu_0)= \max_{\oo} (u-\varphi)$.
 
  We say that a  lower semicontinuous
 function $w :\oo \mapsto \R$ is a
{\rm{ viscosity subsolution}} of \eqref{e.hj} if for
 every test function $\varphi \in \cC_s(\oo)$
 we have
$$
 -\partial_t \varphi(t_0,\mu_0) 
 +H(\mu_0,\partial_\mu \varphi(t_0,\mu_0)) \ge \ell(t_0,\mu_0),
$$
 at every $(t_0,\mu_0) \in \cO$ satisfying
 $(w-\varphi)(t_0,\mu_0)= \min_{\oo} (w-\varphi)$.
 
A function $v: \oo \mapsto \R$ is a {\rm{viscosity solution}}
 if its lower semicontinuous envelope $v^*$
 is a subsolution, and its lower semicontinous envelope $v_*$ is a
 subsolution.
 \end{definition}

\begin{remark}
\label{rem:2} {\rm{In view of \eqref{eq:c2star}, if $\varphi$ is a 
test function, then $\pmu \varphi(t,\mu) \in \cC^2$ with its 
derivatives satisfying $\pmu \varphi(t,\mu), 
|\nabla \pmu \varphi(t,\mu)|^2 \in \cC_*$, 
and $D^2 \pmu \varphi(t,\mu) \in \cC_b$.
Note that these test functions
are not necessarily bounded and 
may grow quadratically.  As our analysis is in
the Wasserstein space $\cP_2$,
this relaxation is natural,
and is utilized in the comparison
proof.
}}\end{remark}

The following is standard and is proved in \cite{SY}.

\begin{corollary}
\label{c.value}
Under Assumption \ref{a.main},
the dynamic programming holds. Consequently, $v$
is a viscosity solution of  \eqref{e.hj},
and for each $n$, $v_n$ is a viscosity solution of
$$
-\partial_t v(t,\mu) + H(\mu, \partial_\mu v(t,\mu)) = \ell_n(t,\mu),
\qquad \text{on}\ \ (0,T) \times \cP_2.
$$
\end{corollary} 

\section{Comparison}
\label{sec:main} 

Our main result is the comparison for the 
Eikonal equation \eqref{e.hj},
and  its proof is given 
later in this section.  Recall that the state space is
$\oo=[0,T]\times \cP_2(\R^d)$, and we endow
it with the product of Euclidean and weak$^*$ topologies.

\begin{theorem}
\label{t.compare} Suppose that Assumption \ref{a.main} holds,
$u: \oo \mapsto \R$ is an upper
semi-continuous, bounded viscosity sub-solution 
of \eqref{e.hj}, and $w :\oo \mapsto \R$ is a lower
semi-continuous, bounded viscosity super-solution  of \eqref{e.hj}.
Further assume that $u(T,\cdot) \le w(T,\cdot)$.
Then, $u \le w$ on $\oo$.  In particular,
the value function $v$ is the
unique  continuous, bounded  viscosity solution of the 
dynamic programming equation \eqref{e.hj}
and the terminal condition $v(T,\cdot)=g$.
\end{theorem}

We start with a simple computation and estimates.
Recall the test functions 
$\cC_s(\oo)$
of Definition~\ref{def:smooth},
$n_*, \varrho$ of \eqref{eq:star},
and the Fourier basis $e(x,\xi)$.
\begin{lemma}
\label{lem:derivative} 
For $\eta \in \cM(\R^d)$,
set $\psi(\eta): = \frac12 \varrho^2(\eta)$ .
Then, for $\mu, \nu \in \cP_2$, 
$$
\kappa(x):=\partial_\mu \psi(\mu-\nu)(x)
= \int_{\R^d}\  (1+|\xi|^2)^{-n^*}\ \fF(\mu-\nu)(\xi) e(x,\xi)\ \d \xi,
\qquad x \in \R^d.
$$
Moreover, $\|\kappa\|_{n^*}= \varrho(\mu-\nu)$.
\end{lemma}
\begin{proof}
Fix $\mu, \nu \in \cP_2$ and set $\eta=\mu-\nu$.
The explicit from of $\kappa:= \partial_\mu \psi(\eta)$ follows
form  a  straightforward computation.  Then,
$\fF(\kappa)(\xi) = \fF(\eta)(\xi) (1+|\xi|^2)^{-n^*}$,
and $\|\kappa\|_{n^*}= \varrho(\eta)$.

\end{proof}

\noindent
{\emph{Proof of Theorem \ref{t.compare}}.}

We complete the proof 
in several steps.  Recall 
$\vt(\mu)=\mu(q)$ defined in \eqref{eq:q}.
Then, $\vt$ is weak* lower-semicontinuous on $\cP_2$,
and  any sublevel set $\{ \mu \in \cP_2 \ :\ \vt(\mu) \le c\}$
is compact.
\vspace{4pt}

\noindent
\emph{Step 1} (\emph{Set-up}). 
Let $u, w$ be as in the statement of the theorem.
Towards a contraposition
suppose that $\sup_{\oo} (u-w)>0$. 
Let $v$ be the value function.
Then, 
$$
0< \sup_{\oo}(u-w) \le \sup_{\oo}(u-v) +\sup_{\oo}(v-w).
$$
Hence,  either $\sup_{\oo} (u-v)>0$, or $\sup_{\oo} (v-w)>0$,
or both must hold.  We analyze the first case and this analysis
can be followed  \emph{mutatis mutandis} to prove the other case.

For a small constant $\gamma_0$, set $\bu(t,\mu):= u(t,\mu)- 2 \gamma_0(T-t+1)$.
We  first fix $\gamma_0$ satisfying $\sup_{\oo} (\bu -v) >0$.  
We then fix $n$ sufficiently large so that
\begin{equation}
\label{eq:eta}
- \partial_t \bu(t,\mu)+H(\mu,\pmu \bu(t,\mu))
\le \ell(t,\mu)  -2\gamma_0 \le \ell_n(t,\mu) -\gamma_0,
\end{equation}
and  $\bu(T,\cdot) \le g-2 \gamma_0 \le g_n$.
In the remainder of
the proof we fix  $\gamma_0, n$ as above.
Next,  set $l := \sup_{\oo}(\bu-v_n)/3$.
\vspace{4pt}

\noindent
\emph{Step 2} (\emph{Doubling the variables}).
Set $\cX=\oo \times \oo$.
For $\eps, \delta >0$, and  $(t,\mu,s,\nu) \in \cX$,
define
$$
\Psi_{\eps,\delta}(t,\mu,s,\nu):=\bu(t,\mu)-v_n(s,\nu)-\frac{1}{2\eps}
((t-s)^2+\varrho^2(\mu-\nu)) -\delta \vt(\mu)- \eps \vt(\nu).
$$
By the previous step, there is
 $(t_0,\mu_0) \in \oo$ such that
$$
2l \le (\bu-v_n)(t_0,\mu_0)
= \Fed(t_0,\mu_0,t_0,\mu_0) + \delta \vt(\mu_0) +\eps \vt(\mu_0).
$$
Then, for all $0<\eps \le \delta \le \delta_*:=l/(2 \vt(\mu_0) +1)$,
$\max_{ \cX}\Fed\geq l >0$.
In the remainder of this proof, we 
always assume that  $\eps \le \delta \le \delta_*$.

Let $(t_k,\muk, s_k,\nuk)$ be a maximizing sequence of $\Fed$.
Since $\bu, v_m$ are bounded,
$$
\delta \vt(\muk) +\eps \vt(\nuk) \le \
(\|\bu\|_\infty + \|v_n\|_\infty)=:c_*.
$$
As the sub-level sets of $\vt$
are compact, 
the sequences $\muk$, $\nuk$ have limit points. 
Since additionally, $v_m, \varrho$ are  continuous, 
and $\bu, -\ \vt$ are upper-semicontinuous, $\Fed$ is also
upper-semicontinuous, and
these limit points achieve 
the maximum of $\Fed$. Hence,  there exists 
a quadruple $(\ted,\mued, \sed,\nued) \in  \cX$ satisfying,  
$\Fed(\ted,\mued, \sed,\nued)
=\max_{ \cX}
\Fed\geq l >0$.  Set 
$$
\etaed:=\mued-\nued,\qquad
\taued:= \ted-\sed
$$  
Then, we also have
\be
\label{eq:bound}
\frac{1}{2\eps}(\taued^2+\varrho^2(\etaed))
+ \delta \vt(\mued)+\eps \vt(\nued)
\le \bu(\ted,\mued)-v_n(\sed,\nued)
\leq c_*.
\ee

\noindent
\emph{Step 3} (\emph{Norm estimate}).
We now use the Lipschitz estimate \eqref{eq:lip} of $v_n$
to obtain a uniform bound for $\varrho(\etaed)/\eps$.
Note that $n$ is already chosen and remains fixed throughout the proof.
As  $\Fed(\ted,\mued,\sed,\mued) \le \Fed(\ted,\mued,\sed,\nued)$, we have
\begin{align*}
u(\ted,\mued)&- v_n(\sed,\mued)- \frac{1}{2\eps} \taued^2
- \delta \vt(\mued)  -\eps \vt(\mued)\\
\le&  \ u(\ted,\mued)-v_n(\sed,\nued)-\frac{1}{2\eps} (\taued^2+\varrho^2(\etaed))
-  \delta \vt(\mued) -\eps \vt(\nued).
\end{align*}
Then, by Proposition \ref{pr:samuel} and \eqref{eq:bound},
\begin{align*}
\frac{1}{2\eps} \varrho^2(\etaed)& \le
v_n(\sed,\mued) - v_n(\sed,\nued) + \eps (\vt(\mued)-\vt(\nued)) \\
& \le \hat{c}_n \varrho(\etaed) + \eps \vt(\mued)
\le \hat{c}_n \varrho(\etaed) + \eps \frac{c_*}{\delta}.
\end{align*}
Therefore, there is a constant $\hc$ depending
only on $\hc_n, c_*$
such that for all $0<\eps,\delta \le 1$,
\begin{equation}
\label{eq:main1}
\frac{\varrho(\etaed)}{\eps}\le \ \frac{\hc}{\sqrt{\delta}}.
\end{equation}

\noindent
\emph{Step 4} (\emph{Letting $\eps$ to zero}).
By \eqref{eq:bound}, $\vt(\mued)\le c_*/\delta$.
Therefore, for each $\delta \in (0,\delta_*]$
there are subsequences
$\{(\ted,\mued)\}\subset \oo$, 
$\{\ted\}\subset   [0,T]$, 
denoted by $\eps$ again, and  limit points
$(\sd,\mud) \in   \oo$, $\td \in [0,T]$
such that as $\eps\downarrow 0$,
$\mued \rightharpoonup\mud$,
$\ted \rightarrow  \td$, and
$\sed\rightarrow \sd$.
By \eqref{eq:bound}, it is clear that $\td=\sd$, 
and 
$\lim_{\eps \downarrow 0} \ \varrho(\mued-\nued) =0$.
We now use Lemma \ref{lem:equiv} 
to conclude that as $\eps \downarrow 0$,
we also have $\nued \rightharpoonup\mud$.

As $\bar{u}(T,\cdot) \le g_n = v_n(T,\cdot)$, 
if  $\td$ were to be equal to $T$,
we would have
\begin{align*}
0<l &\le \liminf_{\eps \downarrow 0}\,
\Phi_{\eps,\delta}(\ted,\mued,\sed,\nued)
\le \liminf_{\eps \downarrow 0}\,
[\bar{u}(\ted,\mued) -v_n(\sed,\nued)]\\
&\le \bar{u}(T,\mud)- v_n(T,\mud) \le 0.
\end{align*}
Hence, $\td<T$ and consequently, both
$\ted<T$, and $\sed <T$ for all sufficiently
small $\eps>0$.
\vspace{5pt}

\noindent
\emph{Step 5} (\emph{Viscosity property}).
Set 
\begin{align*}
\fed(t,\mu)&:=\frac{1}{2 \eps}
((t-\sed)^2+\varrho^2(\mu-\nued)) +\delta \vt(\mu),\\
\ped(s,\nu)&:=-\, \frac{1}{2 \eps}
((\ted-s)^2+\varrho^2(\mued-\nu)) -\eps \vt(\nu).
\end{align*}
By Lemma~\ref{lem:derivative}, both 
$\partial_\mu\fed(t,\mu), \partial_\mu \ped(t,\mu) \in \cC^2_*$.  
Hence, $\fed, \ped \in \cC_s(\oo)$, i.e., they
are smooth test functions
in the sense of Definition \ref{def:smooth}.
By using Lemma~\ref{lem:derivative},
we calculate that
$$
\partial_\mu \fed(\ted,\mued)=\ked + \delta q,\qquad
\partial_\nu \ped(\sed,\nued)= \ked -\eps q,
$$
where $q$ is as in \eqref{eq:q}, and for $x \in \R^d$,
$$
\ked(x):= \frac{1}{\eps}
\int_{\R^d} \ (1+ |\xi |^2)^{-n_*}\ \fF(\etaed)(\xi) 
\, e(x,\xi)\, \d \xi\quad \Rightarrow\quad
\|\ked\|_{n^*} =\frac{1}{\eps}\varrho(\etaed).
$$
It is clear that, $\bar{u}(t,\mu)- \fed(t,\mu)$
is maximized at $(\ted, \mued)$.  Since $\ted<T$, 
$\fed \in \cC_s(\oo)$ and $\bar{u}$ is a viscosity 
subsolution of \eqref{eq:eta}, 
$$
- \frac{\ted-\sed}{\eps} +
H(\mued,\ked+ \delta q) \le \ell_n(\ted,\mued)-\gamma_0.
$$
By the viscosity property of $v_n$, a similar
argument implies that
$$
- \frac{\ted-\sed}{\eps} +
 H(\nued, \ked- \eps q)\ge \ell_n(\sed,\nued).
$$

\noindent
\emph{Step 6} (\emph{Estimation}).
We subtract the above inequalities 
 to arrive at
\begin{align*}
0 < \gamma_0 \le&
H(\nued, \ked-\eps q)
 - H(\mued,\ked+ \delta q)
 + \ell_n(\ted,\mued)-\ell_n(\sed,\nued)\\
=:&  \ \Ied+\Jed+\cK_{\eps,\delta},
\end{align*}
where
\begin{align*}
\Ied &:= \frac12 (\mued(\Delta (\ked+\delta q) - 
 \nued(\Delta (\ked-\eps q),\\
\Jed&:=  \frac12 ( \nued(|\nabla \ked-\eps q|^2) -
 \mued(|\nabla(\ked +  \delta q)|^2)),\\
\cK_{\eps,\delta}&:= \ell_n(\ted,\mued)-\ell_n(\sed,\nued).
\end{align*}
By Assumption \ref{a.main}, $\cK_{\eps,\delta}$ converges to zero as 
$\eps \downarrow 0$.  
Moreover, since $\Delta q \le d$, for $\eps \le \delta$,
$$
\Ied =-\frac{1}{2\eps}\, \int_{\R^d}\, \frac{|\xi|^2}{(1+ |\xi|^2)^{n_*}}\, 
|\fF(\etaed)(\xi)|^2\,  \d \xi\, -\,  \frac12(\delta \mued+\eps \nued)(\Delta q) \le \delta d.
$$
Hence,
$0< \gamma_0 \le \Jed + \delta d$.
\vspace{5pt}

\noindent
\emph{Step 7} (\emph{Estimation of $\Jed$}).
In view of Lemma \ref{lem:derivative},
\eqref{eq:sobolev}, Lemma \ref{lem:derivative},
and \eqref{eq:main1},
$$
\|\ked\|_{\cC^1(\R^d)} \le k_d \|\ked\|_{n^*} 
=  k_d \frac{\varrho(\etaed)}{\eps}
\le \frac{k_d  \hat{c}}{\sqrt{\delta}}.
$$
Since $\nabla q(x) = x/q(x)$, $|\nabla q| \le1$,
and by algebra,
$$
\nued(|\nabla(\ked -\eps q)|^2-|\nabla\ked |^2 )
= -\nued(\nabla(2 \ked -\eps q)\cdot \eps \nabla q)
\le \eps (\frac{2k_d  \hat{c}}{\sqrt{\delta}}+\eps).
$$
Similarly,
$$
 \mued(|\nabla\ked |^2 -|\nabla(\ked +  \delta q )|^2)
=  - \mued(\nabla(2 \ked +\delta  q)\cdot\delta \nabla q)
\le  2k_d  \hat{c} \sqrt{\delta}.
$$
 Therefore,
\begin{align*}
\Jed & =  \frac12 ( \nued(|\nabla \ked-\eps p_{\eps,\delta}|^2) -
 \mued(|\nabla(\ked +  \delta q+\eps q_{\eps,\delta})|^2))\\
& \le \frac12 \etaed(|\nabla \ked|^2) + \bar{c} ( \sqrt{\delta} + \eps),
\end{align*}
for some constant $\bar{c}$ independent of $\eps$. 

We have shown that as $\eps \downarrow 0$,
$\mued, \nued \ \rightharpoonup\ \mud$.
In particular, $\mued, \nued$ are tight sequences
and $\etaed \rightharpoonup 0$.
Additionally, since $\|\ked\|_{\cC^1(\R^d)}$
is uniformly bounded, on a subsequence $\ked$ is locally uniformly
convergent.  These imply that $\etaed(|\nabla \ked|^2)$
converges to zero as $\eps \downarrow 0$. Therefore,
$$
\liminf_{\eps \downarrow 0}\ \Jed 
\le  \bar{c}  \sqrt{\delta}.
$$

\noindent
\emph{Step 8} (\emph{Conclusion}). By the previous steps, for every $\delta >0$ 
the following holds,
$$
0 <  \gamma_0 \le  \limsup_{\eps \downarrow 0}\ 
\Jed +\delta d \le \bar{c}  \sqrt{\delta} +\delta d .
$$
Since $\gamma_0>0$, we obtain a contradiction 
by letting $\delta \downarrow 0$.
Hence, $\sup_{\oo} (u-w) \le 0$.

\qed
\appendix
\section{Convergence of measures in $\varrho$ }
\label{app:B}

For any $s>d/2$, any finite Borel
measure is an element of the Sobolev space $\cH_{-s}(\R^d)$.
Hence,  $\varrho=\|\cdot\|_{-n^*}$ is a metric on $\cP_2(\R^d)$.
Although $(\cP_2, \varrho)$ is not complete,
convergence in this space is equivalent to the 
weak$^*$ convergence in the following sense.
\begin{lemma}
\label{lem:equiv}
Assume that a sequence of probability measures
$\muk$ converge to a probability measure $\mu$
in the weak$^*$ topology, i.e., $\muk \rightharpoonup \mu$.  
Then, $\lim_k \varrho(\muk-\mu)=0$.
Additionally, if a sequence of probability
measures $\nuk$ satisfies $\lim_k \varrho(\nuk-\muk)=0$,
then $\nuk \rightharpoonup \mu$ as well.
\end{lemma}
\begin{proof}
As $\muk \rightharpoonup \mu$,
$\lim_{\eps \downarrow 0} \fF(\muk-\mu)(\xi)= 0$ for every $\xi$.
Then, dominated convergence implies 
$\lim_{\eps \downarrow 0}\varrho(\muk-\mu) =0$,
and 
$\limsup_{\eps \downarrow 0}\varrho(\nuk-\mu)
\le \lim_{\eps \downarrow 0}\varrho(\nuk-\muk)
+ \lim_{\eps \downarrow 0}\varrho(\muk-\mu) =0$.
Hence, $\lim_k \nuk(f)=\mu(f)$
for every $f \in \cH_{n^*}(\R^d)$. Then, by a direct approximation argument, we conclude that $\lim_k \nuk(f) =\mu(f)$
for every \emph{compactly} supported $f \in \cC_b(\R^d)$. Hence, $\nuk$ converges 
to $\mu$ vaguely.  
To  prove that they also converge  in the weak$^*$ topology,
we first fix a smooth function
$ h :[0,\infty) \mapsto [0,1]$
satisfying $h(r) =1$ for all $r \in [0,1]$,
and $h(r)=0$ for all $r \ge 2$. For $m>1$, set 
$h_m(r):= h(r/m)$.  Then, for any $f \in \cC_b(\R^d)$, and $m>1$,
\begin{align*}
 | \nuk(f) -\mu(f)|  & \le 
| \nuk(fh_m) -\mu(fh_m)|
+ \nuk(|f(1-h_m|))
+ \mu(|f(1-h_m)|)\\
&\le 
| \nuk(fh_m) -\mu(fh_m)| +
 \|f\|_\infty ( \nuk(1-h_m) + \mu(1-h_m)).
\end{align*}
Since $h_m$ is compactly supported, $\lim_k \nuk(fh_m)
= \mu(fh_m)$.  Therefore, for every $m>1$,
$$
\limsup_k | \nuk(f) -\mu(f)|   \le 
\|f\|_\infty\ (\limsup_k |\nuk(1-h_{m})| + \mu(1-h_m)).
$$
Moreover, $\lim_k \nuk(1-h_m)=1- \lim_k \nuk1(h_m)
= 1-\mu(h_m)=\mu(1-h_m)$, and as $1-h_m$ converges to zero pointwise,
$\lim_m \mu(1-h_m)=0$.  
Hence, we conclude that $\lim_k \nuk(f)=\mu(f)$ for every
$f \in\cC_b(\R^d)$, and consequently, $\muk$ converges to $\mu$
in the weak$^*$ topology.

\end{proof}

\section{Proposition \ref{pr:samuel}}

Here, we outline the proof of Proposition \ref{pr:samuel}
in several steps.  We fix $n$ and set
$$
L(t,\mu,x) := \partial_\mu \ell_n(t,\mu)(x),\qquad
G(\mu,x) := \partial_\mu g_n(\mu)(x),
\qquad (t,\mu,x) \in \oo \times \R^d.
$$

\noindent{\emph{Step 1.}} (\emph{Reformulation}).
The optimal 
control problem is in fact a deterministic control
problem which has an equivalent representation.
Indeed, for a given initial condition $(t_0,\mu_0) \in \oo$,
let $\cA(t_0,\mu_0)$ be the set of all pairs $(\alpha,m)$ 
satisfying,
\begin{itemize}
\item $m : [t_0,T] \mapsto \cP_2$ is continuous with $m(t_0,\cdot)=\mu_0$;
\item  $\alpha: [t_0,T]\times \R^d \mapsto \R^d$ is Borel measurable and
$\int |\alpha(t,x)|^2 m(t,\d x)\ \d t <\infty$;
\item for any $\phi \in \cC(\R^d)$,
$$
\int_{\R^d} \phi(x) m(s,\d x) = \mu(\phi)+ \int_{t_0}^s 
\int_{\R^d} \ ( \frac12 \Delta \phi(x) + \alpha(t,x) \cdot \nabla \phi(x))\
m(t,\d x)\ \d t.
$$
\end{itemize}
The final condition simply
states that $m(t,\cdot)$ is the law of 
a solution to the stochastic differential equation
$dX_t = \alpha(t,X_t)\d t + \d W_t$.

Then, the value function has the following
equivalent representation \cite{Dau}(Section 2),
$$
v_n(t_0,\mu_0) =\inf_{(\alpha, m) \in \cA(t_0,\mu_0)}
\ \int_{t_0}^T [\ \ell_n(t,m(t,\cdot)) +\frac12 \int_{\R^d} |\alpha(t,x)|^2
m(t,\d x)\ ] \ \d t + g_n(m(T,\cdot)).
$$

\noindent{\emph{Step 2.}}  (\emph{Smooth optimal feedback control}). 
By Pontryagin maximum principle
(see Theorem 2.2 of \cite{Dau}
with constraint $\Psi \equiv 0$),
for any initial condition $(t_0,\mu_0)$
there exists an optimal pair $(\alpha^*,m^*) \in \cA(t_0,\mu_0)$.
Moreover, $\alpha^*(t,x)= - \nabla u(t,x)$
where $u$ is the solution of the following 
Eikonal equation, 
$$
-\partial_t u(t,x) -\frac12 \Delta u(t,x) 
+\frac 12 |\nabla u(t,x)|^2 = \hat{L}(t,x):= L(t,m^*(t,\cdot),x),
\qquad (t,x) \in (0,T) \times \R^d,
$$
with the final condition $u(T,x)= \hat{G}(x):= G(m^*(T,\cdot),x)$.

Recall that $L,G$ have continuous
and bounded derivatives of order $2n^*$.
By standard elliptic regularity (see Lemma 3.1 \cite{DDJ}),
the solution $u$ of the above equation 
satisfies $u(t,\cdot) \in \cC^{2n^*}_b(\R^d)$ 
with norms uniformly bounded in time. We may then 
rewrite the above equation as
$$
-\partial_t u(t,x) -\frac12 \Delta u(t,x) 
+\frac12 A(t,x) \cdot \nabla u(t,x) = \hat{L}(t,x),
$$
where $A(t,x):= \nabla u(t,x)$.  We 
now know that $A(t,\cdot) \in \cC^{2n^*-1}_b(\R^d)$.
Also by hypothesis $\hat{L}(t,\cdot), \hat{G}$
are in $\cH_{2n^*}(\R^d)$.  As the above
equation is linear with smooth coefficients,
standard techniques imply that 
$u(t,\cdot) \in \cH_{2n^*-1}(\R^d)$ with 
norms uniformly bounded in time.  In particular,
we conclude that there is a feedback optimal control
$\alpha^*$ satisfying the estimate
$$\|\alpha^*(t,\cdot)\|_{\cC^{2n^*-1}(\R^d)}
+ \|\alpha^*(t,\cdot)\|_{\cH_{2n^*-1}(\R^d)}
\le C,
$$
with a constant $C$ depending only on the 
norms of $\hat{L}, \hat{G}$. In particular,
$C$ is independent of the initial condition $(t_0,\mu_0)$.
\vspace{4pt}

\noindent{\emph{Step 3.}}
(\emph{Conclusion}).  We now 
follow  \emph{mutadis mutandis} the proofs of Proposition 3.2
and Lemma 3.3 in \cite{DDJ},
(that proves exactly the same 
result on the torus),
to obtain the Lipschitz estimate \eqref{eq:lip}.
 Alternatively, Section 7 of \cite{SY}
also implies the Lipschitz continuity using the 
smoothness of the optimal feedback control.

\bibliographystyle{abbrvnat}
\bibliography{wasserstein}
\end{document}